\def\draft{n}
\documentclass{amsart}
\usepackage{fullpage,amssymb}
\usepackage{epsfig}

\theoremstyle{plain}

\newtheorem{theorem}{Theorem}
\newtheorem*{theoremnn}{Theorem}
\newtheorem{proposition}{Proposition}[section]
\newtheorem{lemma}[proposition]{Lemma}
\newtheorem{corollary}[proposition]{Corollary}
\newtheorem{claim}[proposition]{Claim}

\theoremstyle{definition}
\newtheorem{definition}[proposition]{Definition}

\theoremstyle{remark}
\newtheorem{example}[proposition]{Example}
\newtheorem{exercise}[proposition]{Exercise}

\newtheorem{remark}[proposition]{Remark}

\def\printname#1{
	\if\draft y
		\smash{\makebox[0pt]{\hspace{-0.5in}
			\raisebox{8pt}{\tt\tiny #1}}}
	\fi }

\newcommand{\psdraw}[2]
           {\begin{array}{c} \hspace{-1.3mm}
	\raisebox{-4pt}{\epsfig{figure=draws/#1.eps,width=#2}}
	\hspace{-1.9mm}\end{array}}

\newlength{\standardunitlength}
\catcode`\@=11
\long\def\@makecaption#1#2{%
       \vskip 10pt

\setbox\@tempboxa\hbox{
         \small\sf{\bfcaptionfont #1. }\ignorespaces #2}%
       \ifdim \wd\@tempboxa >\captionwidth {%
           \rightskip=\@captionmargin\leftskip=\@captionmargin
           \unhbox\@tempboxa\par}%
         \else
           \hbox to\hsize{\hfil\box\@tempboxa\hfil}%
       \fi}
\font\bfcaptionfont=cmssbx10 scaled \magstephalf
\newdimen\@captionmargin\@captionmargin=2\parindent
\newdimen\captionwidth\captionwidth=\hsize
\catcode`\@=12

\def\lbl#1{\label{#1}\printname{#1}}

\def\BZ{\mathbb Z}
\def\BQ{\mathbb Q}
\def\BR{\mathbb R}

\def\A{\mathcal A}

\def\O{\mathcal O}
\def\T{\mathcal T}

\def\T{\mathcal T}

\def\La{\Lambda}
\def\l{\lambda}
\def\Ga{\Gamma}
\def\S{\Sigma}

\def\ga{\gamma}

\def\fti{finite type invariant}

\def\AS{\mathrm{AS}}
\def\IHX{\mathrm{IHX}}

\def\lth{\Lambda_{\Theta}}

\def\la{\langle}
\def\ra{\rangle}

\def\y1x{\underset{y1x}\ast}
\def\ti{\widetilde}

\def\eyes{\operatorname{\psfig{figure=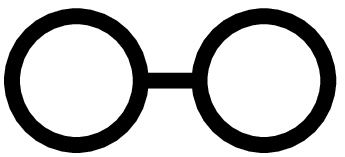,height=0.1in}}}

\def\sminus{\smallsetminus}

\def\deg{\text{deg}}
\def\Zrat{Z^{\mathrm{rat}}}

\def\Lath{\Lambda_{\Theta}}
\def\tiLath{\widetilde{\Lambda}_{\Theta}}

\def\Wh{\mathrm{Wh}}

\def\Aut{\mathrm{Aut}}

\def\Sym{\mathrm{Sym}}

\def\lk{\mathrm{lk}}
\def\b{\beta}
\def\a{\alpha}

\def\e{\epsilon}
\def\mat#1#2#3#4{\left(
\begin{matrix}
    #1 & #2  \\
    #3 & #4
\end{matrix}
\right)}
\def\pt{\partial}
\def\sub{\subset}
\def\hair{\mathrm{Hair}}
\def\longto{\longrightarrow}

\newcommand{\Id}{{1\!\!1}}
\def\down{{^\downarrow}}
\def\up{{^\uparrow}}

\def\strutb#1#2#3{\overset{#1}{\underset{#2}{
\begin{array}{c} \vspace{0.0cm}
\uparrow
\vspace{-0.25cm} \\        
| \vspace{-0.45cm} \\      
\bullet \vspace{0.00cm}   
\end{array} }}\! #3}

\def\lkZ{\mathrm{lk}_{\BZ}}

\begin{document}

\title{On Knots with trivial Alexander polynomial}
\author{Stavros Garoufalidis}
\address{Department of Mathematics \\
               University of Warwick \\
               Coventry, CV4 7AL, UK. }
\email{stavros@maths.warwick.ac.uk}
\author{Peter Teichner}
\address{Department of Mathematics \\
              University of California in San Diego          \\
              9500 Gilman Drive  \\
              La Jolla, CA, 92093-0112, USA.}
\email{teichner@math.ucsd.edu}

\thanks{The  authors are partially supported by NSF grants
             DMS-02-03129 and DMS-00-72775 respectively. The second 
author was also supported by the Max-Planck Gesellschaft.
             This and related preprints can also be obtained at {\tt
http://www.math.gatech.edu/$\sim$stavros} and {\tt
http://math.ucsd.edu/$\sim$teichner}
\newline 1991 {\em Mathematics Classification.} Primary 57N10. Secondary
57M25.\newline {\em Key words and phrases:} Alexander polynomial,
knot, Seifert surface, Kontsevich integral, concordance, slice, clasper.}

\date{ October 7, 2003 \hspace{0.3cm}
First edition: May 31, 2002.}

\begin{abstract}
We use the 2-loop term of the Kontsevich integral to show that
there are (many) knots with trivial Alexander polynomial which don't
have a  Seifert surface whose genus equals the rank of the Seifert form.
This is one of the first applications of the Kontsevich integral to
intrinsically $3$-dimensional questions in topology.

Our examples contradict a lemma of Mike Freedman, and we explain what
went wrong in his argument and why the mistake is irrelevant for topological
knot concordance.
\end{abstract}

\maketitle


\section{A question about classical knots}
\lbl{sec.question}

Our starting point is a wrong lemma of Mike Freedman in
\cite[Lemma~2]{F1}, dating back before his proof of the
$4$-dimensional topological
Poincar\'{e} conjecture. To formulate the question, we  need the following

\begin{definition} \lbl{def.minrank} A knot in 3-space has {\em minimal Seifert
rank} if it has a Seifert surface whose genus equals the rank of the
Seifert form.
\end{definition}
Since the Seifert form minus its transpose gives the
(nonsingular) intersection form on the Seifert surface, it follows
that the genus
is indeed the smallest possible rank of a Seifert form. The formula
which computes
the Alexander polynomial in terms of the Seifert form shows that
knots with minimal
Seifert rank have trivial Alexander polynomial. Freedman's wrong lemma claims
that the converse is also true. However, in the argument he overlooks
the problem
that  S-equivalence does {\em not} preserve the condition of minimal
Seifert rank. It
turns out that not just the argument, but also the statement of the
lemma is wrong.
This has been overlooked for more than 20 years, maybe because none
of the  classical
knot invariants can distinguish the subtle difference between trivial Alexander
polynomial and minimal Seifert rank.

In the last decade, knot theory was overwhelmed by a plethora of new
``quantum''invariants, most notably the HOMFLY polynomial (specializing to
the Alexander and the Jones polynomials), and the Kontsevich integral.
Despite their rich structure, it is not clear how strong these invariants are
for solving open problems in low dimensional topology. It is the
purpose of this
paper to provide one such application.

\begin{theorem}\lbl{thm.main}
There are knots with trivial Alexander polynomial
which don't have minimal Seifert rank. More precisely, the 2-loop part of the
Kontsevich integral induces an  epimorphism
$\overline{Q}$ from the monoid of knots with trivial Alexander
polynomial, onto an
{\em infinitely generated} abelian group, such that
$\overline{Q}$ vanishes on knots with minimal Seifert rank.
\end{theorem}
The easiest counterexample is shown in Figure~\ref{fig.example},
 drawn using surgery on a clasper. Surgery on a clasper is a refined
form of Dehn surgery (along an embedded trivalent graph, rather than
an embedded link) which we explain in Section \ref{sec.claspers}.
Clasper surgery is an elegant way of drawing knots that amplifies the 
important features of our example suppressing irrelevant information
(such as the large number of crossings of the resulting knot).
For example, in Figure~\ref{fig.example}, if one pulls the central edge of the
clasper out of the visible Seifert surface, one obtains an S-equivalence to a
nontrivial knot with minimal Seifert rank.

\begin{figure}[htpb]
$$
\psdraw{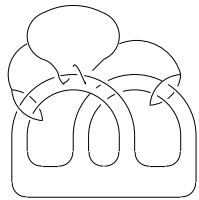}{1in}
$$
\caption{The simplest example, obtained by a clasper surgery on  the
unknot.}\lbl{fig.example}
\end{figure}

\begin{remark} All of the above notions make sense for knots in
homology spheres.
Our proof of Theorem~\ref{thm.main} works in that setting, too.
\end{remark}

Since \cite[Lemma 2]{F1} was the starting point of what eventually
became Freedman's
theorem that all knots with trivial Alexander polynomial are
topologically slice, we
should make sure that the above counterexamples to his lemma don't  cause any
problems in this important theorem. Fortunately, an argument independent of the
wrong lemma can be found in \cite[Thm. 7]{F2}, see also
\cite[11.7B]{FQ}. However, it uses unnecessarily the surgery exact
sequence and some facts from $L$-theory.

In an appendix, we shall give a more direct proof that Alexander
polynomial~1 knots
are topologically slice. We use no machinery, except for a single
application of
Freedman's main disk embedding theorem \cite{F2} in $D^4$. To
satisfy  the
assumptions of this theorem, we employ a triangular base change for
the intersection
form of the complement of a Seifert surface in $D^4$, which works for
all  Alexander
polynomial~1 knots. By Theorem~\ref{thm.main}, this base change does
{\em not}  work
on the level of Seifert forms, as Freedman possibly tried to anticipate.

\section{A relevant quantum invariant}
\lbl{sec.relevant}

The typical list of knot invariants that might find its way into a text book or
survey talk on {\em classical} knot theory, would  contain the
Alexander polynomial,
(twisted) signatures, (twisted) Arf invariants, and maybe knot determinants. It
turns out that all of these invariants can be  computed from the
homology of the
infinite cyclic covering of the knot complement. In particular, they
all vanish if
the Alexander polynomial is trivial. This condition also implies that certain
``noncommutative'' knot invariants vanish,  namely all those
calculated from the
homology of {\em solvable} coverings of the knot  complement, like the
Casson-Gordon
invariants \cite{CG} or the von Neumann signatures of
\cite{COT}. In fact, the latter are concordance invariants  and, as
discussed above,
all knots with trivial Alexander polynomial are  topologically slice.

Thus it looks fairly difficult to study knots with trivial Alexander
polynomial
using classical invariants. Nevertheless, there are very natural  topological
questions about such knots like the one explained in the previous
section. We do
not know a classical treatment of that question, so we turn to
quantum invariants.

One might want to use the Jones polynomial, which often distinguishes
knots with
trivial Alexander polynomial. However, it is not clear which knots it
distinguishes,
and which values it realizes, so the Jones polynomial
is of no help
to this problem. Thus, we are looking for a quantum invariant that
relates well to classical topology, has good realization properties,
and is one step beyond the Alexander polynomial.

In a development starting with the Melvin-Morton-Rozansky conjecture
and  going all
the way to the recent work of \cite{GR} and \cite{GK1}, the
Kontsevich integral has
been reorganized in a rational form $\Zrat$ which is closer to the algebraic
topology of knots. It is now a theorem  (a restatement of the MMR
Conjecture) that
the ``1-loop'' part of the  Kontsevich integral  gives the same
information as the
Alexander polynomial \cite{BG, KSA}.

The quantum invariant in Theorem~\ref{thm.main} is the ``2-loop''
part $Q$ of the
{\em rational invariant} $\Zrat$ of
\cite{GK1}. We consider $Q$ as an invariant of Alexander polynomial~1 knots
$K$ in integral homology spheres $M^3$, and summarize its properties:
\begin{itemize}
\item
$Q$ takes values in the abelian group
$$
\lth:= \frac{\BZ[t_1^{\pm 1},t_2^{\pm 1},t_3^{\pm 1}]}{(
t_1t_2t_3-1,\quad\Sym_3\times\Sym_2) }
$$ The second relations are given by the symmetric groups
$\Sym_3$ which acts by permuting the $t_i$, and $\Sym_2$ which
inverts the $t_i$
simultaneously.
\item Under connected sums and  orientation-reversing, $Q$ behaves as follows:
\begin{align*} Q(M\# M',K \# K') &= Q(M,K) + Q(M',K')\\ Q(M,-K) &=
Q(M,K) =-Q(-M,K)
\end{align*}
\item If one applies the augmentation map
$$
\e:\lth\to \BZ, \quad t_i\mapsto 1,
$$ then $Q(M,K)$ is mapped to the Casson invariant $\lambda(M)$, normalized by
$\l(S^3_{\text{Right Trefoil},+1})=1$.
\item
$Q$ has a simple behavior under surgery on
{\em null claspers}, see Section~\ref{sec.Q}.
\end{itemize} All these properties are proven in \cite{GR} and in \cite{GK1}.

\begin{proposition}[Realization]
\lbl{prop.onto} Given a homology sphere $M^3$, the image of $Q$ on
knots in $M$ with
trivial Alexander polynomial is the subspace $\e^{-1}(\lambda(M))$ of $\lth$.
\end{proposition}

\begin{remark}
\lbl{rem.onto} The realization in the above proposition is concrete,
not abstract.
In fact, to realize the subgroup $\e^{-1}(\lambda(M))$ one only needs
(connected
sums of) knots which are obtained as follows:  Pick a standard
Seifert surface $\Sigma$
of genus one for the unknot in $M$, and do a surgery along a clasper
$G$ with one
loop and two leaves which are  meridians to the bands of $\Sigma$, just like in
Figure~\ref{fig.example}.  The loop of $G$ may
intersect $\Sigma$ and these intersection create the interesting
examples. Note that all
of these knots are ribbon which implies unfortunately that the
invariant $Q$ does
{\em not} factor through knot concordance, even though it  vanishes
on knots of the
form $K\# -K$.
\end{remark}

Together with the following finiteness result, the above realization
result proves
Theorem~\ref{thm.main}, even for knots in a fixed homology sphere.

\begin{proposition}[Finiteness]
\lbl{prop.finite} The value of $Q$ on knots with minimal Seifert rank
is the subgroup
of $\lth$, (finitely) generated by the three elements
$$ (t_1-1),\quad (t_1-1)(t_2^{-1}-1), \quad (t_1-1)(t_2-1)(t_3^{-1}-1).
$$ This holds for knots in 3-space, and one only has to add
$\lambda(M)$  to all
three elements to obtain the values of $Q$ for knots in a homology sphere $M$.
\end{proposition}

\begin{corollary}
\lbl{cor.vas} If a knot $K$ in $S^3$ has minimal Seifert rank, then
$Q(S^3,K)$ can
be  computed in terms of three Vassiliev invariants of degree~$3,5,5$.
\end{corollary}

The $Q$ invariant can be in fact calculated on many classes of
examples. One  such computation was done in \cite{Ga}: The (untwisted)
Whitehead double of a  knot $K$
has minimal Seifert rank and $K\mapsto Q(S^3,\Wh(K))$ is a nontrivial
Vassiliev invariant of degree~$2$.

\begin{remark} Note that $K$ has minimal Seifert rank if and only if
it bounds a
certain grope of class~3. More precisely, the bottom surface of this
grope is just
the Seifert surface, and the second stages are embedded disjointly
from the Seifert
surface. However, they are allowed to intersect each other. So this
condition is
quite different to the notion of a ``grope cobordism''  introduced in
\cite{CT}.
\end{remark}

In a forthcoming paper, we will study related questions for boundary
links.  This is
made possible by the rational version of the Kontsevich integral for
such links recently defined in \cite{GK1}. The analogue of knots with trivial
Alexander polynomial are called {\em good boundary links}. In
\cite[11.7C]{FQ} this term was used for boundary links whose free cover
has trivial homology. Unfortunately, the term was also used in \cite{F1} for
a class of
boundary links which should be rather called {\em boundary links of minimal
Seifert rank}. This class of links is relevant because they form the atomic
surgery problems
for topological $4$-manifolds, see Remark~\ref{rem.atomic}. By
Theorem~\ref{thm.main}
the two definitions of good boundary links in the literature actually differ
substantially (even for knots). One way to resolve the
``Schlamassel'' would be to
drop this term all together.

\tableofcontents

\section{S-equivalence in homology spheres}
We briefly recall some basic notions for knots in homology spheres.
We decided to include the proofs because they are short and might not be
well known for homology spheres, but we claim no originality. Let
$K$ be a knot  in a homology sphere $M^3$. By looking at the inverse image of
a regular  value under a map $M\sminus K\to S^1$, whose homotopy class
generates
$$
[M\sminus K,S^1] \cong H^1(M\sminus K;\BZ) \cong H_1(K;\BZ) \cong \BZ
\quad \text{ (Alexander duality in $M$)}
$$
one constructs a {\em Seifert surface} $\Sigma$ for $K$. It is a
connected oriented surface embedded in $M$ with boundary $K$. Note that a
priori the resulting surface is not connected, but one just ignores the
closed components. By the  usual discussion about twistings near $K$, one
sees that a collar of $\Sigma$ always  defines the linking number zero
pushoff of $K$.

To discuss uniqueness of Seifert surfaces, assume that $\Sigma_0$ and
$\Sigma_1$ are both connected oriented surfaces in $M$ with boundary $K$.

\begin{lemma} \lbl{lem.uniqueness}
After a finite sequence of ``additions of tubes'', i.e.\ ambient 0-surgeries,
$\Sigma_0$ and $\Sigma_1$ become isotopic.
\end{lemma}

\begin{proof}
Consider the following closed surface in the product $M \times I$ (where
$I=[0,1]$):
$$ \Sigma_0\cup (K \times I) \cup \Sigma_1 \, \subset M \times  I
$$ As above, relative Alexander duality shows that this surface bounds an
connected oriented
$3$-manifold $W^3$, embedded in
$M \times I$. By general position, we may assume that the projection
$p:M \times I\to I$ restricts to a Morse function on $W$. Moreover, the usual
dimension counts show that after an ambient isotopy of $W$ in $M
\times I$ one can
arrange  for $p:W\to I$ to be an ordered Morse function, in the sense that the
indices of the  critical points appear in the same order as their
values under $p$.
This can be done  relative to $K
\times I\subset W$ since $p$ has no critical points there.

Consider a regular value $a\in I$ for $p$ between the index 1 and
index 2 critical
points. Then $\Sigma:=p^{-1}(a) \subset M \times \{a\} = M$ is a
Seifert  surface for $K$.
By Morse theory, $\Sigma$ is obtained from $\Sigma_0$ by
\begin{itemize}
\item A finite sequence of small $2$-spheres $S_i$ in $M$ being born,
disjoint from
$\Sigma_0$. These correspond to the index~0 critical points of $p$.
\item A finite sequence of tubes $T_k$, connecting the $S_i$ to (each
other and)
$\Sigma_0$. These correspond to the index~1 critical points of $p$.
\end{itemize} Since $W$ is connected, we know that the resulting
surface $\Sigma$ must be connected.
In case there are no index~0 critical points, it is easy to see that
$\S$ is obtained from $\S_0$ by additions of tubes. We will now reduce
the general case to this case. This reduction is straight forward if the first
tubes $T_i$ that are born have exactly one end on $S_i$, where $i$
runs through all index~0 critical points. Then a sequence of
applications of the {\em lamp cord trick} (in other words, a sequence
of Morse cancellations) would show that up to isotopy one can ignore
these pairs of critical points, which include all index~0 critical points.

To deal with the general case, consider the level just after all
$S_i$ were born and add ``artificial'' thin tubes (in the complement
of the expected $T_k$) to obtain a connected surface. By the lamp cord trick,
this surface is isotopic  to $\Sigma_0$, and the $T_k$ are now tubes on
$\Sigma_0$, producing a connected surface $\Sigma_0'$. Since by
construction the tubes $T_k$ do not go through the artificial tubes,
we can cut the artificial tubes to move from $\Sigma_0'$ back to
$\Sigma$ (through index~2 critical points).

We can treat $\Sigma_1$ exactly as above, by turning the Morse
function upside down, replacing index~3 by index~0, and index~2 by index~1
critical points. The result is a surface $\Sigma_1'$, obtained from
$\Sigma_1$ by adding tubes,  and such that $\Sigma$ is
obtained from $\Sigma_1'$ by cutting other tubes.

Collecting the above information, we now have an ambient Morse
function with only critical points of index~1 and 2, connecting $\Sigma_0$
and $\Sigma_1$ (rel $K$), and a middle surface $\S$ which is tube equivalent
to $\S_0$ and $\S_1$. The result follows.
\end{proof}

The above proof motivates the definition of
S-equivalence, which is  the
algebraic analogue, on the level of Seifert forms, of the geometric addition of
tubes. Given a Seifert surface $\Sigma$ for $K$ in $M$, one defines
the {\em Seifert form}
$$ S_\Sigma: H_1\Sigma \times H_1\Sigma \to \BZ
$$ by the formula $S_\Sigma(a,b):=\lk(a,b\down)$. These are the usual
linking  numbers for
circles in
$M$ and $b\down$ is the circle $b$ on $\Sigma$, pushed slighly off
the  Seifert surface
(in a direction given by the orientations). The downarrow reminds us
that  in the
case of
$a$ and $b$ being the short and long curve on a tube, we are pushing
$b$ {\em into} the tube, and hence the resulting linking number is one.

It should be clear what it means to ``add a tube'' to the Seifert
form $S_\Sigma$: The
homology increases by two free generators $s$ and $l$ (for ``short''
and ``long''
curve on the tube), and the linking numbers behave as follows:
$$
\lk(s,s\down)=\lk(l,l\down)=
\lk(l,s\down)=\lk(s,a\down)=0,\quad\lk(s,l\down)=1,
\quad\forall a\in H_1\Sigma.
$$ Note that there is no restriction on the linking numbers of $l$
with  curves on
$\Sigma$, reflecting the fact that the tube can wind around $\Sigma$
in an arbitrary way.

Observing that isotopy of Seifert surfaces gives isomorphisms of  their Seifert
forms, we are lead to the following algebraic notion. It abstracts
the  necessary
equivalence relation on Seifert forms coming from the non-uniqueness
of the  Seifert
surface.

\begin{definition}\lbl{def.S} Two Seifert surfaces (for possibly
distinct knots) are
called {\em  S-equivalent} if their Seifert forms become isomorphic
after a finite
sequence of  (algebraic) additions of tubes.
\end{definition}

\section{Geometric basis for Seifert surfaces}
It is convenient to discuss Seifert
forms in terms of their  corresponding matrices. So for a given basis
of $H_1\Sigma$,
denote by $SM_\Sigma$ the matrix of linking numbers describing the
Seifert form $S_\Sigma$. For
example, the addition of a tube has the following effect on a Seifert
matrix $SM$:
$$ SM \mapsto \left(
\begin{matrix} SM & 0 & \rho \\
      0 & 0 & 1 \\
      \rho^T & 0 & 0\\
\end{matrix}
\right)
$$ Here we have used the short and long curves on the tube as the
last  two basis
vectors (in that order). $\rho$ is the column of linking number of
the long curve
with the basis elements of $H_1\Sigma$ and $\rho^T$ is its transposed
row. It is clear
that in general this operation can destroy the condition of having
minimal Seifert
rank as  defined in Definition~\ref{def.minrank}. An important invariant of
S-equivalence is the {\em Alexander  polynomial}, defined by
\begin{equation}
\lbl{eq.delta}
\Delta_K(t):=\det( t^{1/2}\cdot SM-t^{-1/2} SM^T)
\end{equation} for any Seifert matrix $SM$ for $K$. One can check that this is
unchanged under S-equivalence, it lies in $\BZ[t^{\pm 1}]$ and satisfies the
symmetry relations $\Delta_K(t^{-1})=\Delta_K(t)$ and
$\Delta_K(1)=1$.

\begin{definition}\lbl{def.basis} Let $\Sigma$ be a Seifert surface
of genus $g$. The
following basis of
$H_1\Sigma$ will be useful.
\begin{itemize}
\item A {\em geometric basis} is a set of embedded simple closed curves
$\{s_1,\dots,s_g, \ell_1,\dots,\ell_g\}$ on $\Sigma$ with the following
geometric  intersections
$$ s_i\cap s_j =\emptyset= \ell_i\cap \ell_j, \text{ and } s_i\cap
\ell_j =\delta_{i,j}
$$ Note that the Seifert matrix $SM_\Sigma$ for a geometric basis
always satisfies
$$ SM_\Sigma - SM_\Sigma^T = \mat{0}{\Id}{-\Id}{0}.
$$
\item A {\em trivial Alexander basis} is a geometric basis such that the
corresponding Seifert matrix can be written in terms of four blocks
of $g \times
g$-matrices as follows:
$$
\mat{0}{\Id+U}{U^T}{V}
$$ Here $U$ is an upper triangular matrix (with zeros on and below
the  diagonal),
$U^T$ is its transpose, and $V$ is a symmetric matrix with zeros
on the diagonal.
\item A {\em minimal Seifert} basis is a trivial Alexander basis such that the
matrices $U$ and $V$ are zero, so the Seifert matrix looks as
simply as could be:
$$
\mat{0}{\Id}{0}{0}
$$
\end{itemize}
\end{definition}

By starting with a disk, and then adding tubes according to the
matrices $U$ and $V$,
it is clear that  any matrix for a trivial Alexander basis can occur
as the Seifert
matrix for the unknot. The curves $s_i$ above are the short curves on
the tubes, and
$\ell_j$ are  the long curves. The matrix $U$ must be lower
triangular because the long
curves can  only link those short curves that are already present.
The following
lemma explains our choice of notation above:

\begin{lemma}
\lbl{lem.basis} Any Seifert surface has a geometric basis. Moreover,
\begin{itemize}
\item A knot has trivial Alexander polynomial if and only if there is  Seifert
surface with a trivial Alexander basis.
\item A knot has minimal Seifert rank if and only if it has a
Seifert surface with
a minimal Seifert basis.
\end{itemize}
\end{lemma}

\begin{proof} By the classification of surfaces, they always have a geometric
basis. If a knot has a trivial Alexander basis, then an elementary
computation using
Equation \eqref{eq.delta} implies that it has trivial Alexander
polynomial. Finally,
the Seifert matrix for a minimal Seifert basis obviously has minimal rank.

So we are left with showing the two converses of the statements in
our lemma. Start
with a knot with trivial Alexander polynomial. Then by Trotter's
theorem \cite{Tr}
it is S-equivalent to the unknot, and hence its Seifert form is
obtained  from the
empty form by a sequence of algebraic additions of tubes. Then an
easy  induction
implies that the resulting Seifert matrix $SM_\Sigma$ is as claimed,
so we are  left with
showing that the corresponding basis can be chosen to be geometric on
$\Sigma$.  But
since $SM_\Sigma - SM_\Sigma^T$ is the standard (hyperbolic) form, we
get a symplectic
isomorphism of $H_1\Sigma$ which sends the given basis into a
standard (geometric) one.
Since  the mapping class group realizes any such symplectic
isomorphism, we see that
the given  basis can be realized by a geometric basis.

Finally, consider a Seifert surface with minimal Seifert rank. By
assumption, there
is a basis of $H_1\Sigma$ so that the Seifert matrix looks like
$$
SM_\Sigma=\mat{0}{A}{0}{B}
$$
Since $\Delta(1)=1$, Equation \eqref{eq.delta} implies that $A$ must be
invertible, and hence there is a base change so that the Seifert matrix has the
desired form
$$
SM_\Sigma=\mat{0}{\Id}{0}{0}
$$
Just as above one shows that this matrix is also realized by a geometric basis.
\end{proof}

\begin{corollary}
\lbl{cor.Matveev} Every knot in $S^3$ with minimal Seifert rank $g$ can be
constructed from a standard genus $g$ Seifert surface of the unknot,
by tying the
$2g$ bands into a 0-framed string link with trivial linking numbers:
$$
\psdraw{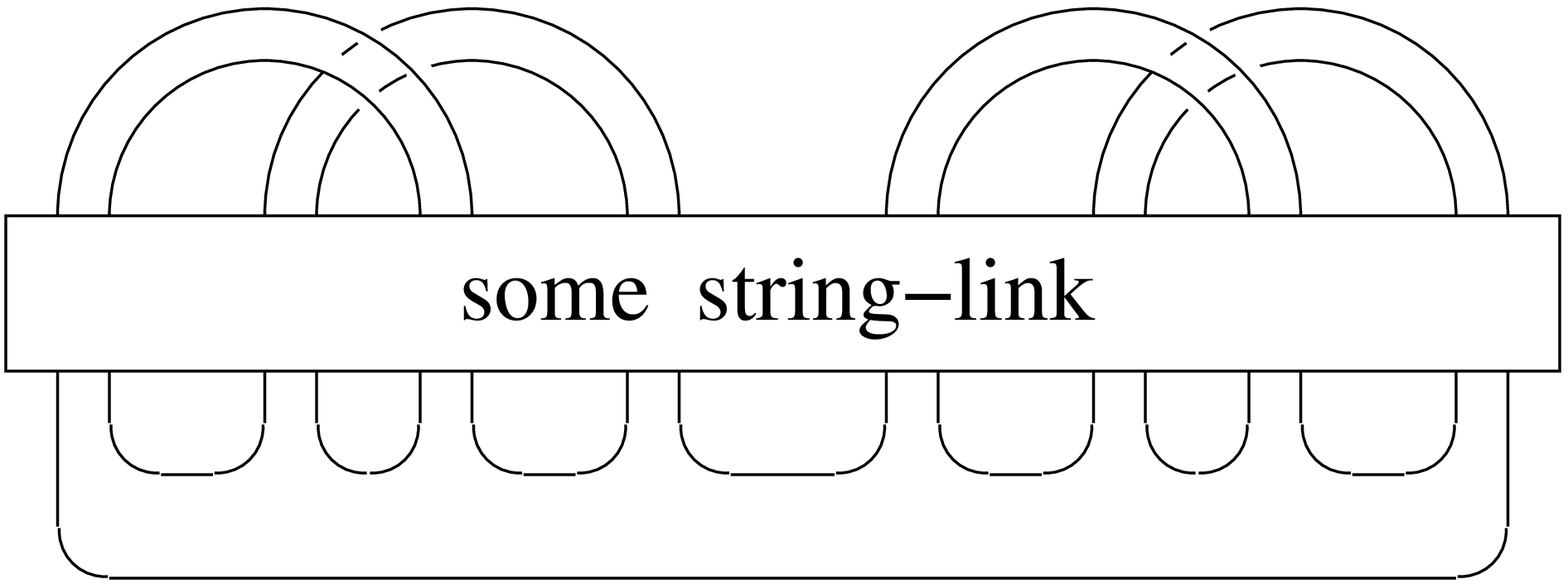}{2.5in}
$$
\end{corollary}

\section{Clasper Surgery}
\lbl{sec.claspers}

As we mentioned in Section \ref{sec.question}, we can construct examples
of knots that satisfy Theorem \ref{thm.main} using {\em surgery on claspers}.
Since claspers play a key role in geometric constructions, as well as
in realization of quantum invariants, we include a brief discussion here.
For a reference on claspers\footnote{By clasper we mean precisely the
object called {\em clover} in \cite{GGP}. For the sake of Peace in the World, 
after the Kyoto agreement of September 2001 at RIMS, we  decided to follow 
this terminology.}  and their
associated surgery, we refer  the reader to \cite{Gu2,H} and also
\cite{CT,GGP}.

{\em Surgery} is an operation of cutting, twisting and pasting within the
category of smooth manifolds. A low dimensional example of surgery is the
well-known {\em Dehn surgery}, where we start from a framed link $L$ in
a 3-manifold $M$, we cut out a tubular neighborhood of $L$, twist the boundary
using the framing, and glue back. The result is a 3-dimensional manifold
$M_L$.

Clasper surgery is entirely analogous to Dehn surgery, excpet that it
is operated on claspers rather than links. A
a clasper is a thickening of a trivalent
graph, and it has a preferred set of loops, called the leaves. The degree of
a clasper is the number of trivalent vertices (excluding those at the
leaves). With our conventions, the smallest clasper is a
Y-clasper (which has degree one and three leaves), so we explicitly
exclude struts (which would be of degree zero with two leaves).

A clasper of degree~1 is an embedding $G: N \to M$ of a regular
neighborhood $N$ of the graph $\Ga$ (with 4 trivalent vertices and 6 edges)

$$
\psdraw{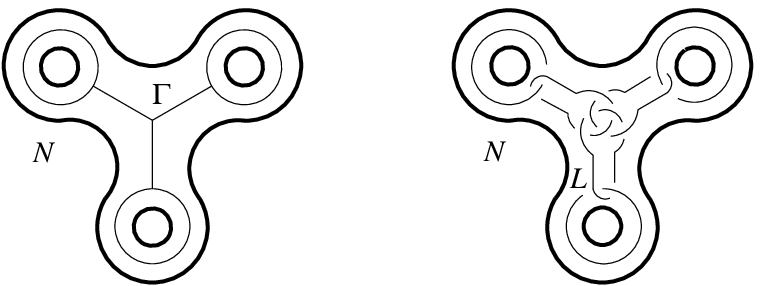}{3in}
$$
into a 3-manifold $M$. Surgery on $G$ can be described by removing the genus~3
handlebody $G(N)$ from
$M$, and regluing by a certain diffeomorphism of its
boundary (which acts trivially on the homology  of the boundary).
We will denote the result of surgery by $M_G$.
To explain the regluing diffeomorphism, we describe surgery on $G$ by
surgery on
the following framed six component link $L$ in $M$: $L$ consists of a
$0$-framed
Borromean ring and an arbitrarily framed three component
link, the so-called {\em leaves} of $G$, see the figure above. The
framings of the
leaves reflect the prescribed neighborhood $G(N)$ of $\Ga$ in $M$.

If one of the leaves is $0$-framed and bounds an embedded disk disjoint
from the rest of $G$, then surgery on $G$ does not change the 3-manifold
$M$, because the gluing diffeomorphism extends to $G(N)$. In terms of the
surgery on $L$ this is explained by a sequence of Kirby
moves from $L$ to the empty link (giving a diffeomorphism $M_G \cong
M$). However,
if a second link $L'$ in $M \smallsetminus G(N)$ intersects the disk bounding
the 0-framed leaf of $L$ then the pairs $(M,L')$ and $(M_G,L')$ might not be
diffeomorphic. This is the way how claspers act on knots or links in a fixed
$3$-manifold $M$, a point of view which is most relevant to this paper.

A particular case of surgery on a clasper of degree~1 (sometimes called a
Y-move) looks locally as follows:

$$
\psdraw{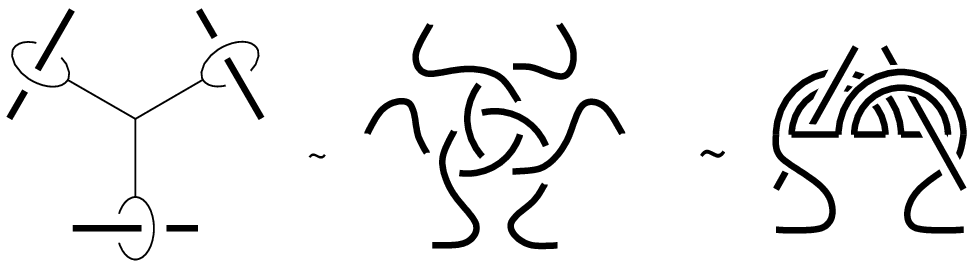}{3.5in}
$$
In general, surgery on a clasper $G$ of degree $n$ is defined in terms
of simoultaneous surgery on $n$ claspers $G_1, \dots, G_n$ of
degree~1. The $G_i$
are obtained from $G$ by breaking its edges and inserting 0-framed Hopf linked
leaves as follows:

$$
\psdraw{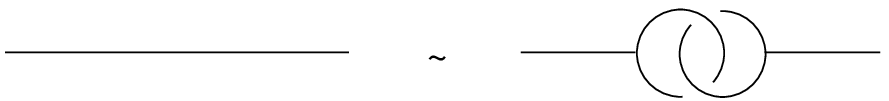}{2in}
$$
In particular, consider the clasper $G$ of degree~2 in Figure
\ref{fig.example}, which has two leaves and two edges. We can insert two
pairs of Hopf links in the edges of $G$ to form two claspers $G_1$ and $G_2$
of degree~1, and describe the resulting clasper surgery on $G_1$ and $G_2$ by
using twice the above figure on each of the leaves of $G$.

\begin{exercise}
\lbl{ex.clasper}
Draw the knot which is described by surgery on a clasper of degree~2
in Figure \ref{fig.example}.
\end{exercise}

It should be clear from the drawing why it is easier to describe knots
by clasper surgery on the unknot, rather than by drawing them explicitly.
Moreover, as we will see shortly, quantum invariants behave
well under clasper surgery.

\section{The $Q$ invariant}\lbl{sec.Q}

\subsection{A brief review of the $\Zrat$ invariant}\lbl{sub.review}

The quantum invariant we want to use for Theorem~\ref{thm.main} is
the Euler-degree $2$ part of the rational invariant $\Zrat$ of
\cite{GK1}. In this section we will give a brief review of the full $\Zrat$
invariant. Hopefully, this will underline the general ideas more clearly, and
will be a useful link with our forthcoming work. $\Zrat$ is a rather
complicated object; however it simplifies when evaluated on Alexander
polynomial $1$ knots, as was explained in \cite[Remark 1.6]{GK1}. In
particular, it is a map of monoids (taking connected sum to multiplication)
$$
\Zrat: \text{Alexander polynomial 1 knots} \longto \A(\La)
$$
where the range is a new algebra of diagrams with beads defined as follows. We
abbreviate the ring of Laurent polynomials in $t$ as $\La:=\BZ[t^{\pm 1}]$.

\begin{definition}
\lbl{def.Ala}
$\A(\La)$ is the completed $\BQ$-vector space generated by pairs
$(G,c)$, where $G$ is a trivalent graph, with oriented edges and  vertices and
$c:\mathrm{Edges}(G)\to \La$ is a $\La$-coloring of $G$,  modulo the
relations:
$\AS$, $\IHX$, Orientation Reversal, Linearity, Holonomy and Graph
Automorphisms, see Figure \ref{relations4} below.
$\A(\La)$ is graded by the {\em Euler degree} (that is, the number of vertices
of graphs) and the completion is
with respect to this grading. $\A(\La)$ is a commutative algebra with
multiplication given by the disjoint union of graphs.
\end{definition}

\begin{figure}[htpb]
$$
\psdraw{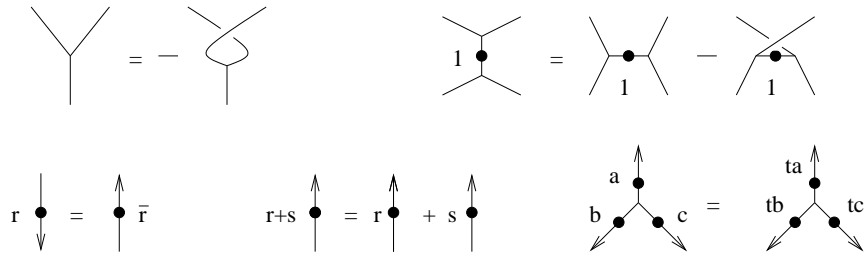}{4.5in}
$$
\caption{The $\AS$, $\IHX$,
Orientation Reversal, Linearity, and Holonomy Relations.}\lbl{relations4}
\end{figure}

 Notice that a connected trivalent graph $G$ has $2n$ vertices, $3n$ edges,
and its Euler degree equals to $-2\chi(G)$, where $\chi(G)$ is the {\em Euler
characteristic} of $G$. This explains the name ``Euler degree''.

Where is the $\Zrat$ invariant coming from? There is an important
{\em hair} map
$$
\hair: \A(\La) \longto \A(\ast)
$$ which is defined by replacing a bead $t$ by an exponential of hair:
$$
\strutb{}{}{t} \mapsto \sum_{n=0}^\infty \frac{1}{n!}
\, \,
\psdraw{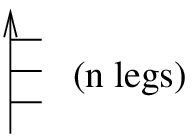}{.7in}
$$
Here, $\A(\ast)$ is the completed (with respect to the {\em Vassiliev
degree}, that is half the number of vertices) $\BQ$-vector
space spanned by vertex-oriented unitrivalent graphs, modulo the $\AS$
and $\IHX$ relations. It was shown in \cite{GK1} that when evaluated on knots
of Alexander polynomial~1, the Kontsevich integral $Z$ is determined by the
rational invariant $\Zrat$ by:
\begin{equation}
\lbl{eq.Zhair}
     Z=\hair \circ \Zrat
\end{equation}
Thus, in some sense $\Zrat$ is a {\em rational lift} of the Kontsevich
integral. Note that although the Hair map above is not 1-1 \cite{P},
the invariants
$Z$ and $\Zrat$ might still contain the same information.
The existence of the $\Zrat$ invariant was predicted by Rozansky, \cite{R},
who constructed a rational lift of the colored Jones function, i.e., for the
image of the Kontsevich integral on the level of the $\mathfrak{sl}_2$ Lie
algebras. The $\Zrat$ invariant was constructed in \cite{GK1}.

How can one compute the $\Zrat$ invariant (and therefore, also the Kontsevich
integral) on knots with trivial Alexander polynomial?  This is a
difficult question; however $\Zrat$ is a graded object, and in each degree it
is a {\em \fti\ } in an appropriate sense.
In order to explain this, we need to recall the {\em null move} of \cite{GR},
which is defined in terms of surgery on a special type of clasper.
Consider a knot $K$ in a homology sphere $M$ and a clasper
$G \sub M\sminus K$ whose leaves are null homologous knots in the knot
complement $X=M\sminus K$.  We will call such claspers {\em null} and will
denote the result of the corresponding surgery by $(M,K)_G$.
Surgery on null claspers preserves the set of Alexander polynomial~1 knots.
Moreover, by results of \cite{Ma} and \cite{MN} one can untie every
Alexander polynomial~1 knot via surgery on some null clasper, see \cite[Lemma
1.3]{GR}.

As usual in the world of \fti s, if $G=\{G_1,\dots,G_n\}$ is a
collection of null
claspers, we set
$$ [(M,K),G]:=
\sum_{I \subset \{0,1\}^n} (-1)^{|I|} (M,K)_{G_I}
$$
where $|I|$ denotes the number of elements of $I$ and $(M,K)_{G_I}$ stands  for
the result of simultaneous surgery on $G_i$ for all $i \in I$.
A {\em \fti\ of  null-type~$k$} by definition vanishes on all such
alternating sums with
$
k < \deg(G):=\sum_{i=1}^n \deg(G_i).
$

\begin{theorem} (\cite{GK1})
\lbl{thm.GKfti1}
$\Zrat_{2n}$ is a \fti\ of null-type $2n$.
\end{theorem}

Furthermore, the degree~$2n$ term (or {\em symbol}) of $\Zrat_{2n}$
can be computed in terms of the  equivariant linking numbers of the leaves
of $G$, as we explain next. Fix an Alexander polynomial~1 knot $(M,K)$,
and consider a {\em null homologous link}
$C \sub X$ of two ordered components, where $X=M\sminus K$.
The lift $\ti C$ of $C$ to the $\BZ$-cover $\ti X$ of $X$ is a link. Since
$H_1(\ti X)=0$ (due to our assumption that $\Delta(M,K)=1$) and
$H_2(\ti X)=0$ (true for $\BZ$-covers of knot complements) it makes  sense to
consider the linking number of $\ti C$. Fix a choice  of lifts $\ti
C_i$ for the components of $C$. The equivariant linking
number is the finite sum
$$
\lkZ(C_1,C_2)=\sum_{n \in \BZ}  \lk(\ti C_1, t^n \, \ti C_2) \, t^n
\, \in\BZ[t^{\pm 1}]=\La .
$$
Shifting the lifts $\ti C_i$ by $n_i\in\BZ$ multiplies this expression by
$t^{n_1-n_2}$. There is a way to fix this ambiguity by considering an
arc-basing of $C$, that is a choice of disjoint embedded arcs $\ga$ in
$M\sminus (K \cup C)$ from a base point to each of the components of $C$.
In that case, we can choose a lift of $C \cup \ga$ to $\ti X$ and define
the equivariant linking number $\lkZ(C_1,C_2)$. The result is independent
of the lift of $C \cup \ga$, but of course depends on the arc-basing $\ga$.

It will be useful for computations to describe an alternative way of fixing
the ambiguity in the definition of equivariant linking numbers. Given $(M,K)$
consider a Seifert surface $\Sigma$ for $(M,K)$, and a link $C$ of two
ordered components in $M\sminus \S$. We will call such links
{\em $\Sigma$-null}. Notice that a $\S$-null link is $(M,K)$-null, and
conversely, every $(M,K)$-null link is $\S$-null for some Seifert surface
$\S$ of $(M,K)$.
Given a $\S$-null link $C$ of two ordered components, one can construct the
$\BZ$-cover $\ti X$ by cutting $X$ along $\Sigma$, and then putting
$\BZ$ copies
of this  fundamental domain together to obtain $\ti X$. It is then obvious that
there  are canonical lifts of $\Sigma$-null links which lie in one fundamental
domain and using them, one can define the equivariant linking number of
$C$ without ambiguity.

This definition of equivariant linking number agrees with the previous one
if we choose basing arcs which are disjoint from $\S$.
\begin{example}
\lbl{ex.levine}
Consider a standard Seifert surface $\Sigma$ for the unknot $\O$. Let
$C_i$ be two meridians of the bands of $\Sigma$; thus $(C_1,C_2)$ is
$\Sigma$-null. If these bands are not dual, then $(\O,C_1,C_2)$ is an unlink
and hence $\lkZ(C_1,C_2)=0$. If the bands are dual, then this 3-component
link is the Borromean rings. Recall that the Borromean rings are the Hopf
link with one component Bing doubled (and the other one being $\O$). Then
one can pull apart that link, in the complement of $\O$, by introducing
two intersections (of opposite sign) between $C_1$ and $C_2$, differing by the
meridian $t$ to $\O$. This shows that in this case
$$
\lkZ(C_1,C_2)= t-1
$$
\end{example}

In order to give a formula for the symbol of $\Zrat_{2n}$, we need to recall
the useful notion of a complete contraction of an $(M,K)$-null clasper
$G$ of degree $2n$, \cite[Sec.3]{GR}. Let $G^{break}=\{G_1,\dots, G_{2n}\}$
denote the collection of degree $1$ claspers $G_i$ which are obtained
by inserting a Hopf link in the edges of $G$.
Choose arcs from a fixed base point to the trivalent vertex of each
$G^{nl}_i$, which allows us to define the equivariant linking numbers
of the leaves of $G^{break}$. Let
$G^{nl}=\{G^{nl}_1,\dots, G^{nl}_{2n}\}$ denote the collection of
abstract unitrivalent graph obtained by removing the leaves of the
$G_i$ (and leaving one leg, or univalent vertex, for each leave behind).
    Then the {\em complete contraction}
$\la G\ra\in\A(\La)$ of $G$ is defined to be the sum over all ways of gluing
pairwise the legs of $G^{nl}$, with the resulting edges of each summand
labelled by elements of $\La$ as follows:
pick orientations of the edges of $G^{nl}$ such that pairs of
legs that are glued are oriented consistently. If two legs $l$ and $l'$ are
glued, with the orientation giving the order, then we attach the bead
$\lkZ(l,l')$ on the edge created by the gluing.

The result of a complete contraction of a null clasper $G$ is a
well-defined element of $\A(\La)$. Changing the edge orientations is taken
care of by the symmetry of the equivariant linking number as well as the
orientation reversal relations. Changing the arcs is taken care by the
holonomy relations in $\A(\La)$.

Then the complete contraction $\la G\ra\in\A(\La)$ of a single clasper
$G$ with $\Sigma$-null leaves is easily checked to be the sum over
all ways of gluing
pairwise the legs of
$G^{nl}$, with the resulting edges of each summand labelled by
elements of $\La$ as
follows: First pick orientations of the edges of $G^{nl}$ such that
pairs of legs
that are glued are oriented consistently. If two legs $l$ and $l'$
are glued, with the
orientation giving the order, then we attach the bead $\lkZ(l,l')$ on the edge
created by the gluing. In addition, each internal edge $e$ of
$G^{nl}$ is labelled by
$t^n$, where $n\in\BZ$ is the intersection number of $e$ with the
Seifert surface
$\Sigma$.

One can check directly that this way of calculating a complete
contraction of a clasper
$G$ with $\Sigma$-null leaves is a well-defined element of $\A(\La)$:
Changing the
edge orientations is taken care of by the symmetry of the equivariant
linking number as
well as the orientation reversal relations. The holonomy relations in $\A(\La)$
correspond beautifully to Figure~\ref{fig.surfacehol} in which a
trivalent vertex
of $G$ is pushed through $\Sigma$.

\begin{figure}[htpb]
$$
\psdraw{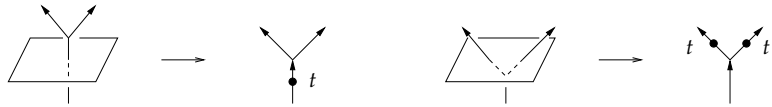}{4in}
$$
\caption{A surface isotopy that explains the Holonomy
Relation.}\lbl{fig.surfacehol}
\end{figure}

Finally, we can state the main result on calculating the invariant $\Zrat$.
\begin{theorem}(\cite[Thm.4]{GK1})
\lbl{thm.GKfti2}
If $(M,K)$ is a knot with trivial Alexander polynomial and $G$ is a
collection of $(M,K)$-null claspers of degree $2n$, then
$$
\Zrat_{2n}([(M,K),G])= \la G \ra \in \A_{2n}(\La)
$$
\end{theorem}

\subsection{A review of the $Q$ invariant}\lbl{sub.reviewQ}

We will be interested in $Q=\Zrat_2$, the loop-degree~2 part of
$\Zrat$. It turns out that $Q$ takes values in a lattice
$\A_{2,\BZ}(\La)$, that is the abelian subgroup of $\A(2,\La)$
generated by integer multiples of graphs with beads. The next lemma
(taken from \cite[Lemma 5.9]{GK2}) explains the definition of $\Lath$.

\begin{lemma}
\lbl{lem.Qvalues} There is an isomorphism of abelian groups:
\begin{equation}
\lbl{eq.la}
\Lath \longrightarrow \A_{2,\BZ}(\La) \hspace{0.5cm} \text{given by:}
\hspace{0.5cm}
\a_1 \, \a_2 \, \a_3 \mapsto
\psdraw{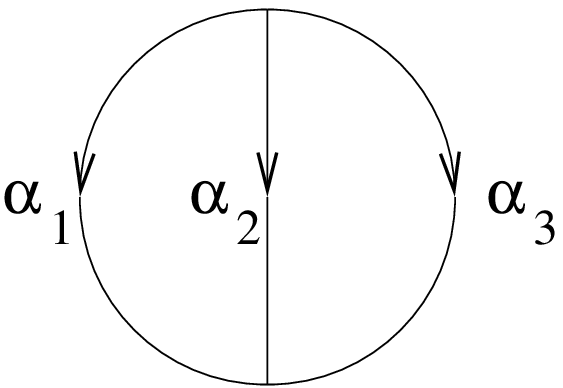}{0.8in}.
\end{equation}
\end{lemma}

\begin{proof} Since $\Aut(\Theta) \cong \Sym_3 \times \Sym_2$, it is easy to
see that the above map is well-defined. There are two trivalent graphs of
degree $2$, namely $\Theta$ and
$\eyes$. Using the Holonomy Relation, we can assume that the labeling
of the  middle
edge of $\eyes$ is $1$. In that case, the $\IHX$ relation implies that
$$
\psdraw{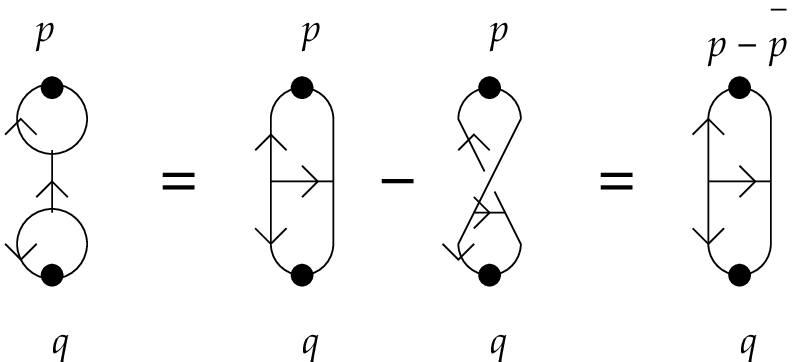}{2in}
$$

This shows that the map in question is onto. It is also easy to see
that it is a
monomorphism.
\end{proof}

Let us define the {\em reduced} groups
$$
\ti\A(\La)=\mathrm{Ker}(\A(\La)\to\A(\phi))
$$ induced by the augmentation map $\e:\La\to\BZ$. Let
$\tiLath:=\mathrm{Ker}(\e: \Lath \to \BZ)$. The proof of the above
lemma implies
that there is an isomorphism:
$$
\tiLath \cong \ti\A_{2,\BZ}(\La).
$$

\subsection{Realization and finiteness}\lbl{sub.realf}

\begin{proof}[Proof of Proposition \ref{prop.onto} (Realization)]
Let us first assume that the ambient 3-manifold $M=S^3$. It is easy to see
that $\ti \Lath$ is generated by $(t_1-1)t_2^n t_3^m$  for $n,m \in \BZ$, so
we only need to realize these values. Consider a standard genus one Seifert
surface $\Sigma$ of an unknot with bands $\{\a,\b \}$ and the clasper $G$
$$
\psdraw{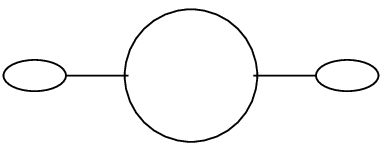}{0.5in}
$$
of degree $2$ (with two leaves shown as ellipses above). Choose an
embedding of $G$
into $S^3 \smallsetminus \O$ in such a way that the two leaves are
$0$-framed meridians of the two bands of $\Sigma$ and the two internal edges of
$G$ intersect $\Sigma$ algebraically $n$ respectively $m$ times. Then
$G$ is a $\Sigma$-null clasper and Theorem \ref{thm.GKfti2}, together with
Exercise~\ref{ex.levine} we get
$$
Q(S^3,\O_G)=-Q([(S^3,\O),G])= (1-t_1)t_2^n t_3^m \in \tiLath.
$$
The realization result follows for $M=S^3$. For the case of a general
homology sphere $M$, use the behavior of $Q$ under connected sums.
To show that the constructed knots are ribbon, we refer to
\cite[Lem.2.1, Thm.5]{GL}, or \cite[Thm.4]{CT}.
\end{proof}

The next lemma gives a clasper construction of all minimal Seifert rank knots.
We first introduce a useful definition. Consider a surface
$\Sigma \subset S^3$ and a clasper $G \subset S^3\sminus \pt \Sigma$.
We say that $G$  is {\em
$\Sigma$-simple} if the leaves of $G$ are $0$-framed  meridians of the
bands of $\Sigma$ and the edges of $G$ are disjoint from $\Sigma$.

\begin{lemma}
\lbl{lem.Yconstruct}
Every knot in $S^3$ with minimal Seifert rank can be
constructed from a standard Seifert surface $\Sigma$ of the unknot,
by surgery on a disjoint collection of $\Sigma$-simple Y-claspers.
\end{lemma}

\begin{proof} The result follows by Lemma~\ref{cor.Matveev} and the
fact, proven by
Murakami-Nakanishi \cite{MN}, that every string-link with trivial
linking numbers
can be untied by a sequence of Borromean moves. In terms of $\O$,
these Borromean
moves are  $\Sigma$-simple Y-clasper surgeries (with the leaves being
$0$-framed meridians
to the bands of $\Sigma$).
\end{proof}

\noindent
{\em Proof of Proposition~\ref{prop.finite}.} (Finiteness)
Consider a knot $K$ in $S^3$ with minimal Seifert rank. By
Lemma~\ref{lem.Yconstruct} it is obtained from a standard Seifert
surface $\Sigma$ of an
unknot $\O$ by surgery on a disjoint collection $G$ of
$\Sigma$-simple Y-claspers. The fact
that
$Q$ is an invariant of type $2$ implies that
$$
Q(S^3,K)=-Q((S^3,\O)-(S^3,\O)_G)=-\sum_{G' \sub G} Q([(S^3,\O),G'])
+\sum_{G'' \sub G} Q([(S^3,\O),G''])
$$
where the summation is over all claspers $G'$ and $G''$ of degree $1$ and
$2$ respectively. The $Q([(S^3,\O),G''])$ terms can be computed by complete
contractions and using Example~\ref{ex.levine}, it follows that they
contribute only summands of the form $(t_i-1)$.

Next we simplify the remaining terms, which are given by
$\Sigma$-simple Y-claspers $G'\sub G$. Note that we can work modulo
$\Sigma$-simple claspers of degree~$>1$ by the above argument. Using the
Sliding Lemma (\cite[Lem.2.5]{GR}) we can move around all edges and finally
put $G'$ into a standard position as in Figure~\ref{fig.remain} below.

\begin{figure}[htpb]
\lbl{fig.remain}
$$
\psdraw{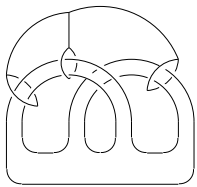}{1in}
$$
\caption{The remaining knots, possibly with half-twists (not shown) on the
edges of the clasper.}\lbl{remain}
\end{figure}

We are reduced to $\Sigma$ of genus one because if the 3 leaves of
$G'$ are meridians to 3 distinct bands of $\Sigma$, the unknot $\O$ would
slip off the clasper altogether, i.e., surgery on the simplified $G'$ does
not alter $\O$.

This means that we are left with a family of 4 examples, given by the various
possibilities of the half-twists in the 3 edges of the the clasper in
Figure~\ref{fig.remain}. Let $\a$ and $\b$ denote the two bands of the
standard genus~1 surface $\S$, and let $m_{\a}, m_{\b}$ (resp.
$\ell_{\a}, \ell_{\b}$) denote the knots which are meridians (resp. longitudes)
of the bands.

Let $G'$ denote the $\S$-simple clasper of degree~1 as in Figure~\ref{remain}.
It has 3 leaves $m_{\a}, m_{\a}$ and $\ell_{\b}$.

\begin{claim}
We have
$$
[(S^3,\O),G']=[(S^3,\O),G'']+[(S^3,\O),G''']
$$
modulo terms of degree~2, where $G''$ is a $\S$-simple clasper with leaves
$m_{\a}, m_{\a}, \ell_{\a}$ and $G'''$ is obtained from $G''$ by replacing
the edge of $\ell_{\a}$ by one that intersects $\S$ once.
\end{claim}

\begin{proof}(of the claim)
Observe that $m_{\b}$ is isotopic to $\ell_{\a}$ by an isotopy rel $\S$.
Use this isotopy to move the leaf $\ell_{\b}$ of $G'$ near the $\a$ handle,
and use the Cutting a Leaf lemma (\cite[Lem.2.4]{GR}) to conclude the proof.
\end{proof}

Going back to the proof of Proposition \ref{prop.finite}, we may apply
the Cutting a Leaf lemma once again to replace $G''$ by a $\S$-simple
clasper with leaves two copies of $m_{\a}$ together with a meridian of
one copy of $m_{\a}$. For this clasper, the surface $\S$ can slide off,
and as a result surgery gives back the unknot. Work similarly for $G'''$,
and conclude that $Q([(S^3,\O),G'])$ lies in the subgroup of $\Lath$
which is generated by the elements
$$
(t^{\e_1}_1-1),\quad (t^{\e_1}_1-1)(t_2^{\e_2}-1), \quad (t^{\e_1}_1-1)
(t^{\e_2}_2-1)(t_3^{\e_3}-1)
$$
for all $\e_i = \pm 1$. Using the relations in $\Lath$, it is easy to
show that this subgroup is generated by the three elements as claimed
in Proposition \ref{prop.finite}. This concludes the proposition for
knots in $S^3$.

In the case of a knot $K$ with minimal Seifert rank in a general
homology sphere $M$, we may  untie it by surgery on a collection of
$\Sigma$-simple Y-claspers, $\Sigma$ a standard Seifert surface for the unknot
$\O$. That is, we may assume that $(M,K)=(S^3,\O)_G$ for some $\Sigma$-null
clasper $G$ whose leaves are meridians of the bands of $\S$ and have
framing $0$ or $\pm 1$. We can follow the previous proof to conclude our
result.
\qed

\begin{proof}[Proof of Corollary \ref{cor.vas}]
As we discussed previously, the
rational invariant $\Zrat$ determines the Kontsevich integral via Equation
\eqref{eq.Zhair}. It follows that
$\hair \circ Q$ is a power series of Vassiliev invariants.  Although
the $\hair$ map
is not 1-1, it is for diagrams with two loops, thus $\hair \circ Q$
determines $Q$.

Consider the image of
$t_1-1$, $(t_1-1)(t_2^{-1}-1)$ and $(t_1-1)(t_2-1)(t_3^{-1}-1)$
under the $\hair$
map in $\A(\ast)$. It follows that the  Vassiliev invariants of
degree $3,5$ and $5$
which separate the uni-trivalent graphs
$$
\psdraw{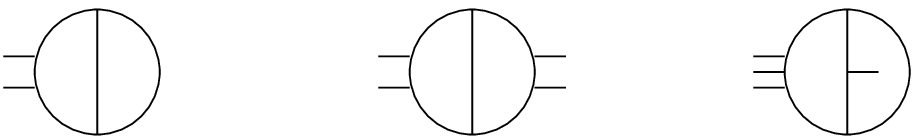}{2in}
$$
determine the value of $Q$ on knots with minimal Seifert rank.
\end{proof}

\appendix

\section{Knots with trivial Alexander polynomial are  topologically slice}

A complete argument for this fact can be found in \cite[Thm.7]{F2}, see also
\cite[11.7B]{FQ}). However, that argument uses unnecessarily the surgery exact
sequence for the trivial as well as infinite cyclic fundamental 
group. Moreover, one
needs to know Wall's surgery groups $L_i(\BZ[\BZ])$ for $i=4,5$.

We shall give a
direct argument in the spirit of
\cite{F1} but without assuming that the knot has minimal Seifert rank 
(which Freedman
did assume indirectly). The simple new ingredient is the triangular 
base change,
Lemma~\ref{lem.triangular}.
    Note that at the time of writing \cite{F1}, the topological disk
embedding theorem was not known, so the outcome of the constructions
below was much weaker than an actual topological slice.

The direct argument uses a single application  of Freedman's
main disk embedding theorem \cite{F2}. In \cite{F2} it is not stated 
in its most
general form which we need here, so we really use the disk embedding theorem
\cite[5.1B]{FQ}. So let's
first recall this basic theorem. It works in any $4$-manifold with {\em good}
fundamental group, an assumption which up to day is not known to be 
really necessary.
In any case, cyclic
groups are known to be good which is all we need in this  appendix.
Note that the
second assumption, on dual $2$-spheres, is well known to be
necessary.  Without this
assumption, the proof below would imply that every  ``algebraically
slice'' knot,
i.e., a knot whose Seifert form has a Lagrangian, is topologically slice. This
contradicts for example the invariants of \cite{CG}. A more direct
reason that this
assumption is necessary was recently given in
\cite{ST}: In the absence of dual $2$-spheres, there are nontrivial secondary
invariants (in two copies of the group ring modulo certain
relations), which are
obstructions to a disk being homotopic to an embedding.

\begin{theoremnn}(Disk embedding theorem \cite[5.1B]{FQ})
Let $\Delta_j:(D^2,S^1)\to (N^4,\partial N)$ be continuous maps of  disks
which are embeddings on the boundary, and assume
that all intersection and  self-intersection numbers vanish in
$\BZ[\pi_1N]$. If $\pi_1N$ is good and there exist algebraically dual
$2$-spheres, then there is a regular homotopy (rel.\ boundary) which takes the
$\Delta_j$ to disjoint (topologically flat) embeddings.
\end{theoremnn}

The assumption on dual $2$-spheres (which is an algebraic condition)
means that there are framed immersions $f_i:S^2\to N$ such that  the
intersection numbers in $\BZ[\pi_1N]$ satisfy
$$
\lambda(f_i,\Delta_j)=\delta_{i,j}
$$ The following simple observation turns out to be crucial for  Alexander
polynomial~1 knots.
\begin{lemma}
\lbl{lem.triangular} There exist algebraically dual $2$-spheres for
$\Delta_i$ if and only if
there exist framed  immersions $g_i:S^2\to N$ with
$$
\lambda(g_i,\Delta_i)=1 \text{ and } \lambda(g_i,\Delta_j)=0 \text{ for } i>j.
$$
So the matrix of intersection numbers of $g_i$ and $\Delta_j$ needs
to have zeros
only below the diagonal.
\end{lemma}
\begin{proof} Define $f_1:=g_1$, and then inductively
$$ f_i:=g_i - \sum_{k<i} \lambda(g_i,\Delta_k) f_k.
$$
Then one easily checks that $\lambda(f_i,\Delta_j)=\delta_{i,j}$.
\end{proof}

\begin{remark} The disk embedding theorem is proven by an application
of another
embedding theorem
\cite[5.1A]{FQ}, to the Whitney disks pairing the intersections among  the
$\Delta_i$. Thus \cite[Theorem 5.1A]{FQ} might be considered as more basic. It
sounds  very similar to \cite[Theorem 5.1B]{FQ}, except that the assumptions on
trivial intersection and self-intersection numbers is moved from the
$\Delta_i$ to
the dual $2$-spheres. Hence one looses  the information about the
regular homotopy
class of $\Delta_i$.

In most applications, one wants this homotopy information, hence we have stated
theorem 5.1B as the basic disk embedding theorem. However, in the
application below we might as well have used 5.1A directly, by
interchanging the
roles of $s_i$ and $\ell_i$.
\end{remark}

The following proof will be given for knots (and slices)  in $(D^4,S^3)$ but it
works just as  well in $(C^4,M^3)$ where
$M$ is any homology sphere and $C$ is {\em the} contractible  topological
$4$-manifold with boundary $M$.

\begin{proof}[Proof of the appendix title] Since the knot $K$ has
trivial Alexander polynomial,  Lemma~\ref{lem.basis} shows that we can choose a
Seifert  surface $\Sigma_1$ with a trivial Alexander basis
$\{s_1,\dots,s_g,\ell_1,\dots,\ell_g\}$. Pick generically immersed disks
$\Delta(s_j)$ (respectively $\Delta(\ell_j)$) in $D^4$ which bound $s_j\down$
(respectively $\ell_j$). So these disks are disjoint on the boundary, and the
intersection numbers satisfy
$$
\Delta(s_i)\cdot\Delta(s_j)=\lk(s_i\down,s_j\down)=\lk(s_i\down,s_j)=0
\quad\text{ and }\quad \Delta(s_i)\cdot\Delta(\ell_j)=\lk(s_i\down,\ell_j).
$$ By Definition~\ref{def.basis}, the latter is a triangular matrix,
which will
turn out to be the crucial fact.

Now we ``push'' the Seifert surface $\Sigma_1$ slightly into $D^4$ to obtain a
surface $\Sigma \subset D^4$, and call
$N$ the complement of (an open neighborhood of) $\Sigma$ in $D^4$.
The basic idea of the proof is to use the disk embedding theorem in
$N$ to show that
$\Sigma$ can be ambiently surgered into a disk which will be a slice
disk for our knot
$K$.

To understand the $4$-manifold $N$ better, note that by Alexander duality
$$
H_1N \cong H^2(\Sigma,\pt \Sigma) \cong\BZ \text{ and } H_2N \cong
H^1(\Sigma,\pt \Sigma)
\cong\BZ^{2g}.
$$
Moreover, a Morse function on $N$ is given by restricting the radius
function on $D^4$. Reading from the center of $D^4$ outward, this
Morse function has
one critical point of index 0, one of index 1 (the minimum of $\Sigma$), and
$2g$ critical points of index 2, one for each band of $\Sigma$.
Together with the above
homology information, this implies that $N$ is homotopy equivalent to
a wedge of a circle
and $2g$ $2$-spheres.

To make the construction of $N$ more precise, we prefer to add an
exterior collar
$(S^3 \times [1,1.5], K\times [1,1.5])$ to $D^4$, i.e. we work with
the knot $K$ in the
4-disk $D_{1.5}$ of radius~$1.5$. Then the pushed in Seifert
surface $\Sigma\subset D_{1.5}$ is just $(K\times [1,1.5])\cup
\Sigma_1$. The normal
bundle of
$\Sigma_1$ in $D_{1.5}$ can then be canonically decomposed as
$$
\nu(\Sigma_1,D_{1.5}) \cong  \nu(S^3,D_{1.5}) \times \nu(\Sigma_1,S^3) =:
\BR_x \times \BR_y
$$
Since $N^4$ is the complement of an open thickening of $\Sigma$
in $D_{1.5}$, we may assume that for points on $\Sigma_1$ the normal
coordinates $x$
vary in the open interval
$(0.9,1.1)$,  and $y$ in $(-\epsilon,\epsilon)$. Here $\epsilon>0$ is
normalized so
that for a curve
$\alpha= \alpha \times 1 \times 0$ on $\Sigma_1$ one has
$$
\alpha \times 1 \times -\epsilon = \alpha\down \quad \text { and }
\quad \alpha \times 1 \times \epsilon = \alpha\up.
$$ Note that by construction, the disks $\Delta(s_j)$ lie in $N$ and
have their
boundary
$s_i\down$ in $\partial N$ and hence one can attempt to apply the
disk embedding
theorem these disks. If we can do this successfully, then the
$\Delta(s_j)$ may be replaced by disjoint embeddings and hence we can surger
$\Sigma$ into a slice disk for our knot $K$.

Let's check the assumptions in the disk embedding theorem:
As mentioned above, $\pi_1N\cong \BZ$ is a good group. By
construction, the (self-)
intersections among the $\Delta(s_j)$ vanish algebraically, even in
the group ring
$\BZ[\pi_1N]$,  because these disks lie in a simply connected part of $N$.

Finally, we need to check that the $\Delta(s_j)$ have algebraically  dual
$2$-spheres. Note that this must be the place where the assumption on
the  Alexander
polynomial is really used, since so far we have only used that
$K$ is ``algebraically slice''. We start with $2$-dimensional tori
$\T_i$ which are the boundaries of small normal bundles of $\Sigma$ in
$D_{1.5}$, restricted to the curves $\ell_i$ in our trivial Alexander
basis of $\Sigma_1$.
More precisely,
$$
\T_i:=\ell_i \times S^1_t \text{ where }  S^1_t:=[0.8,1.2] \times
\{-2 \epsilon,2
\epsilon\} \cup
\{0.8,1.2\} \times [-2 \epsilon,2 \epsilon]
$$
in our normal coordinates introduced above. Note that $S^1_t$ is a
(square shaped)
meridian to $\Sigma$ and freely generates $\pi_1N$.
      By construction, these $\T_i$ lie in our
$4$-manifold $N$. Moreover, they are disjointly embedded and dual to
$\Delta(s_j)$ in the sense that the {\em geometric} intersections are
$$
      \T_i\cap \Delta(s_j) =  (\ell_i\cap s_j) \times (0.8 \times - \epsilon)  =
\delta_{i,j}.
$$ Hence the $\T_i$ satisfy all properties of dual $2$-spheres,
except  that they
are not
$2$-spheres! However, we can use our disks $\Delta(\ell_i)$ with
boundary $\ell_i$ as
follows. First remove  collars $\ell_i \times (0.8,1]$ from these
disks (without
changing their name) so that $\Delta(\ell_i)$ have boundary equal to the
``long curve''
$\ell_i\times 0.8$ on
$\T_i$.  Using two parallel copies of $\Delta(\ell_i)$ we can surger the
$\T_i$ into $2$-spheres $g_i$. These are framed because of our
assumption that the
$\ell_i$ are ``untwisted'', i.e.\ that $\lk(\ell_i,\ell_i\down)=0$
(which is used only
modulo~2). The equivariant intersection numbers are
$$
\lambda(g_i,\Delta(s_j))=\delta_{i,j} +
\Delta(\ell_i)\cdot\Delta(s_j)(1-t)=\delta_{i,j} +\lk(\ell_i,s_j\down)(1-t)
\quad
\in\quad\BZ[\pi_1N]=\BZ[t^{\pm 1}]
$$ because the single intersection point of $\Delta(s_i)$ with $\T_i$
remains and
any geometric intersection point between $\Delta(\ell_i)$ and
$\Delta(s_j)$  is now
turned into exactly two (oppositely oriented) intersections of
$g_i$ with $\Delta(s_j)$. These differ by the group element $t$ going
around the
short curve $S^1_t$ of $\T_i$. By our assumption on the linking  numbers, the
resulting $2$-spheres $g_i$ satisfy the triangular condition from
Lemma~\ref{lem.triangular} and can hence be turned into dual spheres for
$\Delta(s_j)$.

Thus we have checked all assumptions in the disk embedding theorem,
and hence we
may indeed surger $\Sigma$ to a slice disk for $K$ as planned.
\end{proof}

\begin{remark}\lbl{rem.atomic}
Recall that the topological surgery and s-cobordism theorems in
dimension~4 (for all fundamental groups) are equivalent to certain 
``atomic'' links
being free slice
\cite[Ch. 12]{FQ}. These atomic links are all boundary links with
minimal Seifert rank
in the appropriate sense. In particular, if the disk embedding
theorem above was true
for free fundamental groups, then the proof above (without needing
our triangular base
change) would show how to find free slices for all the atomic links.
This shows how one
reduces the whole theory to the disk embedding theorem for free
fundamental groups.
\end{remark}

\ifx\undefined\bysame
	\newcommand{\bysame}{\leavevmode\hbox to3em{\hrulefill}\,}
\fi

\end{document}